\documentclass[12pt]{article}
\usepackage{amsmath,amssymb,amsthm,amscd}

\title{Modularity on vertex operator algebras arising from semisimple
primary vectors}

\author{
  Hiroshi Yamauchi
  \vsb\\
  {\small \textit{Graduate School of Mathematics,}}
  \\
  {\small \textit{University of Tsukuba, Ibaraki 305-8571, Japan}}
  \\
  {\small e-mail: {\sf hirocci@math.tsukuba.ac.jp}}  
  \vsb
}

\date{}

\makeatletter
\@addtoreset{equation}{section}

\theoremstyle{plain}
\newtheorem{thm}{Theorem}[section]
\newtheorem{prop}[thm]{Proposition}
\newtheorem{lem}[thm]{Lemma}
\newtheorem{cor}[thm]{Corollary}

\newtheorem{ass}{Assertion}
\newtheorem{claim}{Claim}
\newtheorem{thm_intro}{Theorem}
\newtheorem{cor_intro}{Corollary}
\newtheorem{lem_intro}{Lemma}

\theoremstyle{definition}
\newtheorem{df}[thm]{Definition}

\theoremstyle{remark}
\newtheorem{rem}[thm]{Remark}

\newcommand{\pf}{\gs{Proof:}\q }

\newcommand{\gs}[1]{\textbf{#1}}
\newcommand{\ds}{\displaystyle}
\newcommand{\fr}[2]{\frac{#1}{#2}}
\newcommand{\dfr}[2]{\dfrac{#1}{#2}}%ams
\newcommand{\cd}{\cdot}
\newcommand{\cds}{\cdots}
\newcommand{\dsum}{\displaystyle \sum}

\renewcommand{\l}{\left}
\renewcommand{\r}{\right}

\newcommand{\vsb}{\vspace{2mm}}

\newcommand{\q}{\quad}
\newcommand{\qq}{\qquad}

\newcommand{\maru}[1]{{\ooalign{\hfil#1\/\hfil\crcr
\raise.167ex\hbox{\mathhexbox20D}}}}

\newcommand{\ruby}[2]{%
 \leavevmode
 \setbox0=\hbox{#1}%
 \setbox1=\hbox{\tiny #2}%
 \ifdim\wd0>\wd1 \dimen0=\wd0 \else \dimen0=\wd1 \fi
 \hbox{%
   \kanjiskip=0pt plus 2fil
   \xkanjiskip=0pt plus 2fil
   \vbox{%
     \hbox to \dimen0{%
       \tiny \hfil#2\hfil}%
     \nointerlineskip
     \hbox to \dimen0{\mathstrut\hfil#1\hfil}}}}

\newcommand{\la}{\langle}
\newcommand{\ra}{\rangle}

\newcommand{\abs}[1]{\lvert{#1}\rvert}

\DeclareMathOperator*{\tensor}{\otimes}

\newcommand{\Z}{\mathbb{Z}}
\newcommand{\C}{\mathbb{C}}

\newcommand{\N}{\mathbb{N}}
\newcommand{\Q}{\mathbb{Q}}

\newcommand{\res}{\mathrm{Res}}
\newcommand{\End}{\mathrm{End}}

\newcommand{\aut}{\mathrm{Aut}}
\newcommand{\wt}{\mathrm{wt}}
\newcommand{\tr}{\mathrm{tr}}
\renewcommand{\hom}{\mathrm{Hom}}
\newcommand{\ch}{\mathrm{ch}}
\newcommand{\SL}{\mathit{SL}}
\newcommand{\id}{\mathrm{id}}

\newcommand{\pii}{\pi i}
\newcommand{\hf}{\fr{1}{2}}

\newcommand{\w}{\omega}
\newcommand{\vacuum}{\mathrm{1\hspace{-3.2pt}l}}
\newcommand{\vac}{\vacuum}

\renewcommand{\pii}{\pi \sqrt{-1}}

\newcommand{\D}{\Delta}
\newcommand{\U}{\mathcal{U}}
\newcommand{\mat}{
  \begin{pmatrix} 
    \alpha & \beta \\ \gamma & \delta  
  \end{pmatrix}
}
\newcommand{\smat}{\big( 
  \begin{smallmatrix} 
    \alpha & \beta \\ 
    \gamma & \delta 
  \end{smallmatrix}
\big)}
\newcommand{\h}{\mathfrak{h}}

\begin{document}

\baselineskip 6mm

\maketitle

\begin{abstract}
  In this article, using an idea of the physics superselection
  principal, we study a modularity on vertex operator algebras
  arising from semisimple primary vectors. 
  We generalizes the theta functions on vertex operator algebras
  and prove that the internal automorphisms do not change 
  the genus one twisted conformal blocks.
\end{abstract}

\section{Introduction} %%%%%%%%%%%%%%%%%%%%%%%%%%%%%%%%%%%%%%%

In the study of the physics superselection principal and its
application in the theory of vertex operator algebras (VOAs), Li
introduced the notion of semisimple primary vectors in \cite{L1}.
Let $V$ be a VOA.  
A vector $u\in V_1$ is called {\it a semisimple primary vector} if it
satisfies: (i) $L(n)u=0$ for $n>0$; (ii) $u_{(m)}u=-\delta_{m,1}\la
u,u\ra \vac$ for $m\geq 0$; (iii) $u_{(0)}$ acts on $V$ semisimply
with rational eigenvalues.
The main feature of a semisimple primary vector $u$ is that it
realizes a functor between distinct categories of irreducible
$V$-modules. 
Define the Delta operator associated to $u$ by
\begin{equation}\label{Delta}
  \Delta (u,z)
  := z^{u_{(0)}}\exp\l( -\dsum_{n=1}^\infty \dfr{u_{(n)}}{n}
  (-z)^{-n}\r) . 
\end{equation}
Let $g$ be an automorphism on $V$ of finite order and denote by
$\sigma (u)$ an {\it internal automorphism} $e^{-2\pii u_{(0)}}$ on $V$.
Assume that $g$ fixes $u$ so that we have $[g,\sigma (u)]=1$.
It is shown in \cite{L2} (see also \cite{L1}) that if
$(W,Y_W(\cd,z))$ is an irreducible $g$-twisted $V$-module, then 
$(W,Y_W(\D (u,z)\,\cd,z))$ is an irreducible $g\sigma (u)$-twisted
$V$-module.
We present a new application of the Delta operators \eqref{Delta} to  
the conformal blocks in the orbifold theory.

Recall the theta functions on VOAs.
They were introduced by Miyamoto in \cite{M1} as a generalization of
the theta functions on lattices.  
Let $u$, $v$ be vectors in weight one subspace of $V$ such that
$u_{(0)}v=0$. 
On every $V$-module $W$, Miyamoto defined the following formal power
series in \cite{M1}:
$$
  Z_W(u;v;\tau )= e^{\pii \la u,v\ra} 
  \tr_W \exp \l( 2\pii  v_{(0)} \r)
  q^{L(0)+u_{(0)}+ \fr{1}{2} \la u,u\ra -\fr{1}{24}c},
$$
where $\la \cd,\cd \ra$ is an invariant bilinear form on $V$ such that
$u_{(1)}v=\la u,v\ra \vac$, $c$ is the central charge of $V$ and $q$
denotes $e^{2\pii\tau}$. 
Based on Zhu's theory \cite{Z}, Miyamoto proved the modularity of
$Z_W(u;v;\tau)$ in \cite{M1}. 
Surprisingly, the modular transformation law of a theta function on a module 
is exactly the same as that of character of the module,  
even though a theta function carries some informations of automorphisms 
and hence those of twisted modules.
On the other hand, using the Delta operators \eqref{Delta}, we can
understand the modularity of $Z_W(u;v;\tau)$ from a different view point.
We know that $(W,Y_W(\Delta (u,z)\, \cd ,z))$ is a $\sigma
(u)$-twisted module.
Denote $(W,Y_W(\D (u,z)\,\cd,z))$ simply by $\tilde{W}$. 
Since the actions of $v_{(0)}$ and $L(0)$ on $\tilde{W}$ can be
identified with those of $v_{(0)}+\la u,v\ra$ and $L(0)+u_{(0)}+ \hf
\la u,u\ra$ on $W$, respectively, the theta functions $Z_W(u;v;\tau)$
can be identified with the trace functions $\tr_{\tilde{W}} \sigma (v) 
q^{L(0)-c/24}$ of an internal automorphism $\sigma (v)$ on a $\sigma
(u)$-twisted module $\tilde{W}$ in the case where both $\sigma (u)$ and 
$\sigma (v)$ are of finite order.
So we can understand a meaning of the modularity of $Z_W(u;v;\tau )$ 
from a view point of orbifold conformal field theory \cite{DLM2}.
Combining results in \cite{M1} and \cite{DLM2}, in this article we will 
reveal the property of the internal automorphisms that they do not change 
the genus one twisted conformal blocks.

Let us explain our results more precisely.
Let $G$ be a finite abelian subgroup of $\aut (V)$.
Assume that $V$ is $C_2$-cofinite and $k$-rational for all $k\in G$. 
Let $u, v$ be mutually commutative rational semisimple
primary vectors (see Sec.{} 4.2) in $V^G$.
For each pair $(g,h)$ in $G\times G$, denote by 
$\{ (W^i(g,h),\phi_i(h)) \mid i=1,\dots,N\}$ the complete set  
of inequivalent irreducible $g$-twisted $h$-stable $V$-modules, 
where $\phi_i (h)$ are (fixed) $h$-stabilizing automorphisms on
$W^i(g,h)$.  
Define the genus one twisted conformal block $\mathcal{C}_1(g,h)$ 
associated to a pair $(g,h)\in G\times G$ as the space of trace functions
$$
  T_{W^i(g,h)} (a,\tau):=\tr_{W^i(g,h)} z^{\wt (a)}Y(a,z) \phi_i (h) 
  q^{L(0)-c/24},\q 1\leq i\leq N.
$$
For $\rho =\smat \in \SL_2(\Z)$, the following modular transformation
is shown in \cite{DLM2}:
\begin{equation}\label{henkei}
  (\gamma \tau +\delta)^{-\wt [a]} T_{W^i(g,h)}(a,\rho \tau)
  = \sum_{j=1}^N A_{ij}(\rho, (g,h)) T_{W^j((g,h)^\rho)}(a,\tau),
\end{equation}
where $A_{ij}(\rho,(g,h))$ are the constants independent of $a$ and
$\tau$, $\rho \tau$ denotes $(\alpha\tau +\beta)(\gamma \tau
+\delta)^{-1}$ and $(g,h)^\rho$ denotes $(g^\alpha h^\gamma,g^\beta
h^\delta)$. 
The matrix $(A_{ij}(\rho,(g,h))_{ij}$ defines a linear isomorphism from
$\mathcal{C}_1(g,h)$ to $\mathcal{C}_1((g,h)^\rho)$.
Define the Schur polynomial $p_s(x_1,x_2,\dots)$ in variables
$x_1,x_2,\dots$, by the equation:
$$
  \exp \l( -\dsum_{n=1}^\infty \dfr{x_n}{n} (-z)^{-n}\r)
  = \dsum_{s=0}^\infty p_s(x_1,x_2,\dots) z^{-s}.
$$
Then define {\it a generalized theta function on $W^i(g,h)$} by
\begin{equation}
  \begin{array}{l}
    \ds
    Z_{W^i(g,h)}(a;(u,v);\tau)
    := \dsum_{s=0}^\infty \tr_{W^i(g,h)}
      \Big\{ (p_s(u_{(1)},u_{(2)},\dots)a)_{(\wt (a)+\lambda_u(a)
       -s-1)}
    \vsb\\
    \qq\qq  \times
      e^{\pii \la u,v\ra} \exp \l( 2\pii v_{(0)} \r) \phi_i(h)^{-1}  
      q^{L(0)+u(0)+\fr{1}{2}\la u,u\ra -c/24}\Big\} ,
  \end{array}
\end{equation}
where $a\in V$, $\lambda_u(a)$ is a scalar such that
$u_{(0)}a=\lambda_u(a)a$ and $q=e^{2\pii\tau}$.
Then our main theorem in this paper is the following:

\begin{thm_intro}
  The generalized theta function $Z_{W^i(g,h)}(a;(u,v);\tau)$
  converges on the upper half plane and gives a vector in
  $\mathcal{C}_1(g\sigma (u),h\sigma (v))$ for each $1\leq i\leq N$.
  Furthermore, we have the following modular transformation for 
  $\rho =\smat\in \SL_2(\Z)$:   
  $$
    \begin{array}{l}
    \ds
     (\gamma \tau +\delta )^{-\wt [a]}
     Z_{W^i(g,h)} \l( a;(u,v);\dfr{\alpha\tau +\beta}{\gamma \tau
      +\delta} \r) 
    \vsb\\
    \ds \q 
    = 
    %(\gamma\tau +\delta )^{\wt [a]}
    \sum_{i=1}^N
    A_{ij}(\rho,(g,h)) \,
    Z_{W^i((g,h)^\rho)}(a;(\alpha u+\gamma v,\beta u+\delta
    v);\tau) ,
  \end{array}
  $$
  where $A_{ij}(\rho,(g,h))$ are the constants given by \eqref{henkei} 
  and independent of $a$, $u$, $v$ and $\tau$.
\end{thm_intro} 

There is an interesting consequence of Theorem 1.
For $\rho\in \SL_2(\Z)$, denote by $\Psi_{(g,h)}(\rho)$ the linear
isomorphism between 
$\mathcal{C}_1(g,h)$ and $\mathcal{C}_1((g,h)^\rho)$ given as
\eqref{henkei}.
Define a linear isomorphism $\Omega_{(g,h)}(u,v): \mathcal{C}_1(g,h)
\to \mathcal{C}_1(g\sigma (u),h\sigma (v))$ by 
$
  T_{W^i(g,h)}(a,\tau) \mapsto 
  Z_{W^i(g,h)}(a;(u,v);\tau) .
$ 
Then we obtain the following corollary.

\begin{cor_intro}
  For each $\rho =
  \smat\in \SL_2(\Z)$, the following diagram commutes: 
  \begin{equation}
  \begin{CD}
    \mathcal{C}_1(g,h)  @>{\ \q \Psi_{(g,h)}(\rho)\q \ }>> 
      \mathcal{C}_1((g,h)^\rho)
    \\
    @V{\Omega_{(g,h)}(u,v)}VV 
    @VV{\Omega_{(g,h)^\rho}(\alpha u+\gamma v,\beta u+\delta v)}V
    \\
    \mathcal{C}_1(g\sigma (u),h\sigma (v)) 
    @>{\ \Psi_{(g\sigma (u),h\sigma (v))}(\rho)\ }>> 
      \mathcal{C}_1((g\sigma (u),h\sigma (v))^\rho) .
  \end{CD}
  \end{equation}
\end{cor_intro}
  
The corollary above says that the $\SL_2(\Z)$-transformation laws for
two different genus one twisted conformal blocks $\mathcal{C}_1(g,h)$ and 
$\mathcal{C}_1(g\sigma (u),h\sigma (v))$ are the same so that 
the internal automorphisms do not changes the genus one twisted conformal
blocks. 
This result seems to suggest that this is true for higher genus conformal
blocks, since internal automorphisms $\sigma (u)$ and $\sigma (v)$ are
generated by the fields in $V$ itself.

The main idea in the proof of Theorem 1 is to use the Delta
operators \eqref{Delta} associated to semisimple primary vectors.
The Delta operator realizes a functor between the category of weak
modules and the category of admissible modules. 
For an admissible $V$-module $(W,Y_W(\cd,z))$, a $V$-module $(W,Y_W(\D 
(u,z)\,\cd,z))$ is not admissible in general.
Namely, even if $V$ is rational, we do not know whether $(W,Y_W(\D
(u,z)\,\cd,z))$ is completely reducible or not. 
Therefore, before we give the proof of Theorem 1, we investigate a
relation between $g$-rationality and $g$-regularity.
We extend the result on the spanning sets for weak modules in
\cite{B} to the twisted modules. 
The following is a simple refinement of Lemma 2.4 of \cite{M3}.

\begin{lem_intro}
  Let $V$ be a $C_2$-cofinite VOA of CFT type
  and $W$ a weak $g$-twisted $V$-module generated by one element $w$.
  Then $W$ is linearly spanned by 
  $$
    \alpha^1_{(n_1)}\cds \alpha^s_{(n_s)}w,\q 
    \alpha^i\in U,\q n_1<\cds <n_s<T,
  $$
  where $U$ is a finite dimensional subspace of $V$ such that
  $V=U+C_2(V)$ and $T$ is a fixed number in $\fr{1}{\abs{g}}\N$.
\end{lem_intro}

It is worth mentioning that the repeat condition in \cite{B} is
now removed by Lemma 1.
As an application, we also extend the result in \cite{ABD} to the
twisted case. 

\begin{cor_intro}
  Every $g$-rational $C_2$-cofinite VOA of CFT type is actually
  $g$-regular.  
\end{cor_intro}

This paper is organized as follows.
In Section 2 we recall basic definitions.
In Section 3 we extend the results in \cite{B} and \cite{ABD} to
the twisted case.
In Section 4.1 we review the theory on the modular invariance
on rational VOAs and in Section 4.2 we review the theory on the
physics superselection principal and semisimple primary vectors. 
Using these theories, we prove Theorem 1 above in Section 4.3.
In Section 4.4 we discuss a relation between our theta functions 
and abelian coset models.

\section{Preliminaries} %%%%%%%%%%%%%%%%%%%%%%%%%%%%%%%%%%%%%%%%%%%%%%

In this paper, we mainly treat VOAs of CFT type.

\begin{df}
  A VOA $V$ is called {\it CFT type} if it has a weight
  decomposition $V=\oplus_{n=0}^\infty V_n$ without negative weights
  and its weight zero subspace is spanned by the vacuum,
  i.e. $V_0=\C \vac$. 
\end{df}

We review the definition of the twisted modules.
Let $g$ be an automorphism on $V$ of finite order $\abs{g}$.
Then we can decompose $V$ as a direct sum of eigenspaces for $g$:
$$
  V=V^0\oplus V^1\oplus\cds \oplus V^{\abs{g}-1},\q \text{where}\ 
  V^r:=\{ a\in V\mid ga=e^{2\pii \, r/\abs{g}}a \} .
$$

\begin{df}
  A weak $g$-twisted $V$-module is a vector space $W$ with a linear
  map 
  $$
    Y_M(\cd,z): a\in V\mapsto Y_M(a,z)=\sum_{n\in\fr{1}{\abs{g}} \Z} 
    a_{(n)} z^{-n-1}\in \End (W)[[z^{\pm\fr{1}{\abs{g}}}]]
  $$
  (called the vertex operator on $M$) satisfying the following:
  \vsb\\
  (i) $Y_M(a,z)=\sum_{n\in \Z} a_{(n+\fr{r}{\abs{g}})}
    z^{-n-1-\fr{r}{\abs{g}}}$ for $a\in V^r$,
  \vsb\\
  (ii) $a_{(n)} w=0$ for $n\gg 0$ where $a\in V$ and $w\in W$,
  \vsb\\
  (iii) $Y_M(\vac,z)=\id_M$,
  \vsb\\
  (iv) The following $g$-twisted Jacobi identity holds for $a\in V^r$
  and $b\in V$: 
  $$
  \begin{array}{l}
    z_0^{-1}\delta\l(\dfr{z_1-z_2}{z_0}\r) Y_M(a,z_1) Y_M(b,z_2)
    %\vsb\\
    %\qq 
      -z_0^{-1}\delta\l(\dfr{-z_2+z_1}{z_0}\r) Y_M(b,z_2) Y_M(a,z_1)
    \vsb\\
    \qq\qq 
    = z_2^{-1}\delta\l(\dfr{z_1-z_0}{z_2}\r) \l(\dfr{z_1-z_0}{z_2}
    \r)^{-\fr{r}{\abs{g}}} Y_M(Y_V(a,z_0)b,z_2) .
  \end{array}
  $$
\end{df}

Let us recall two consequences of the twisted Jacobi identity.
Let $W$ be a weak $g$-twisted $V$-module.
For $a\in V^r$, $b\in V^s$ and $w\in W$, there exists $k\in \N$ such 
that $z^{k+\fr{r}{\abs{g}}} Y_W(a,z)w \in W[[z]]$.
By (iv), we can derive the following associativity (cf.{} \cite{L2}).
\begin{equation}\label{associativity}
  (z_2+z_0)^{k+\fr{r}{\abs{g}}} Y_M(Y_V(a,z_0)b,z_2) w
  = (z_0+z_2)^{k+\fr{r}{\abs{g}}} Y_M(a,z_0+z_2) Y_M(b,z_2) w.
\end{equation}
Let $A$, $B$ be subsets of $V$ and $X$ a subset of $W$.
Set $A\cd X = \la a_{(n)} x\mid a\in A, x\in X, n\in \fr{1}{\abs{g}}
\Z\ra$.
Using \eqref{associativity}, we can show the following
associativity-like relation: 
\begin{equation}\label{associativity-like}
  A\cd (B\cd X) \subseteq (A\cd B)\cd X.
\end{equation}
In particular, $V\cd w$ is a submodule of $W$.

Another consequence of the twisted Jacobi identity is the
iterate formula.
On $V$, there exists $N\gg 0$ such that $(z_1-z_2)^N
Y_V(a,z_1)Y_V(b,z_2) = (z_1-z_2)^N Y_V(b,z_2) Y_V(a,z_1)$.
Then the following formula holds on $W$.
\begin{equation}\label{iterate}
\begin{array}{lll}
  \ds
  (a_{(m)}b)_{(n+\fr{r}{\abs{g}}+\fr{s}{\abs{g}})}
  &=& \ds
    \sum_{i=0}^N \sum_{j=0}^\infty (-1)^j \binom{-r/\abs{g}}{i} 
    \binom{m+i}{j} 
    \vsb\\      
    &&
    \ds \times
    \l\{ a_{(m+i-j+\fr{r}{\abs{g}})}
    b_{(n-i+j+\fr{s}{\abs{g}})} 
    -(-1)^{m+i} b_{(m+n-j+\fr{s}{\abs{g}})} a_{(j+\fr{r}{\abs{g}})}
   \r\} .
\end{array}
\end{equation}

\begin{df}
  An admissible $g$-twisted $V$-module is a weak $g$-twisted
  $V$-module which carries a $\fr{1}{\abs{g}}\N$-grading
  $M=\oplus_{n\in\fr{r}{\abs{g}}\N} M(n)$
  such that $a_{(m)} M(n)\subseteq M(n+\wt (a)-m-1)$ for all $a\in
  V$. 
\end{df}

\begin{df}
  An ordinary $V$-module is a weak $V$-module which carries a
  $\C$-grading $M=\oplus_{s\in \C} M_s$ such that:
  \vsb\\
  (i) $\dim M_s<\infty$,
  \vsb\\
  (ii) $M_{s+N}=0$ for any fixed $s$ and sufficiently small integer 
       $N$, 
  \vsb\\
  (iii) $L(0)w=sw=\wt (w)w$ for $w\in M_s$.
\end{df}

It follows from definitions that every ordinary $g$-twisted $V$-module
is an admissible $g$-twisted $V$-module.
Also, it is shown in \cite{DLM1} that an irreducible admissible
$g$-twisted $V$-module is an irreducible ordinary $g$-twisted
$V$-module. 

\begin{df}
  A VOA $V$ is said to be {\it $g$-rational} if every admissible 
  $g$-twisted $V$-module is a direct sum of irreducible admissible
  $g$-twisted $V$-modules.
  Also, $V$ is said to be {\it $g$-regular} if every weak $g$-twisted
  $V$-module is a direct sum of irreducible ordinary $g$-twisted
  $V$-modules.
  A $1$-rational (resp. $1$-regular) VOA is simply called {\it
  rational} (resp. {\it regular}). 
 
\end{df}

\begin{df}
  For $n\geq 2$, set $C_n(W)= \la a_{(-n)} w \mid a\in V, w\in W\ra$.
  An ordinary $V$-module $W$ is said to be {\it $C_n$-cofinite} if
  $W/C_n(W)$ is of finite dimension. 
\end{df}

There are many conjectures about rationality.
We give some of them below.
\vsb\\
(1) Rationality of orbifold VOA \cite{DVVV}:
if $V$ is rational then $V^G$
is also rational, where $G$ is a finite automorphism group acting on 
$V$ and $V^G$ denotes the fixed point subalgebra 
$\{ a\in V\mid ga=a\ \text{for}\ \text{all}\ g\in G\}$.
\vsb\\
(2) Rationality and $C_2$-cofiniteness \cite{ABD}:
rationality provides $C_2$-cofiniteness.
\vsb\\
(3) Rationality and regularity \cite{DLM3}: 
every rational VOA is regular.
\vsb\\
(4) Relation between rationality and $g$-rationality:
if $V$ is rational, then $V$ is $g$-rational for any finite
automorphism $g$.
\vsb\\
Concerning to the conjecture (1) and (4), we  prove the following.

\begin{prop}\label{g-rationality}
  Let $V$ be a simple VOA and $g$ an automorphism on $V$ of finite
  order $\abs{g}$.
  If the orbifold VOA $V^{\la g\ra}$ is rational, then $V$ is
  $g$-rational.
\end{prop}

\pf
Recall the associative algebras $A_{g,n}(V)$ introduced in \cite{DLM4}.
Since $V^{\la g\ra}$ is rational, all $A_{1,n}(V^{\la g\ra})$, $n \in
\N$, are semisimple by \cite{DLM1}.
Then all $A_{g,n}(V)$, $n\in \fr{1}{\abs{g}}\N$, are also semisimple 
because they are homomorphic images of semisimple algebras 
$A_{1,n}(V^{\la g\ra})$.
Therefore, by a theorem in \cite{DLM4}, $V$ is $g$-rational.
\qed

\begin{rem}
   By the proposition above, if the conjecture (1) is true for arbitrary
   finite cyclic group $G=\la g\ra$, then (4) will follow from (1).
\end{rem}

Recently, Abe, Buhl and Dong proved the following remarkable theorem.

\begin{thm} \label{ABD}
  (\cite[Theorem 4.5]{ABD}) 
  Every rational $C_2$-cofinite VOA of CFT type is
  regular.
\end{thm}
\noindent
This result will be generalized to $g$-twisted case in the next
section.

\section{Spanning set for twisted VOA-modules} %%%%%%%%%%%%%%%

Here we give brief generalizations of the results obtained in \cite{B} 
and \cite{ABD}.

For a VOA $V$ of CFT type, Gaberdiel and Neitzke showed the
following theorem on a spanning set of $V$. 

\begin{thm}\cite{GN}
  Let $V$ be a $C_2$-cofinite VOA of CFT type and write $V=U+C_2(V)$
  with $\dim U<\infty$. 
  Then $V$ is spanned by vectors of the form $\alpha^1_{(-n_1)}\cds
  \alpha^k_{(-n_k)} \vac$, $n_1>\cds >n_k>0$ with each $\alpha^i\in
  U$. 
\end{thm}

We generalize this theorem to weak $g$-twisted $V$-modules.
First, we recall the $g$-twisted universal enveloping algebra
$\U^g(V)$ of $V$ in \cite{DLM1}.
As a tensor product of two vertex algebras $\C
[t^{\pm\fr{1}{\abs{g}}}]$ and $V$, $\hat{V}:= \C [t^{\pm
\fr{1}{\abs{g}}}] \tensor_\C V$ carries a structure of a vertex
algebra and $\mathfrak{g}_V:=\hat{V}/(\fr{d}{dt}\tensor 1+1\tensor
L(-1))\hat{V}$ forms a Lie algebra under the 0-th product induced from 
$\hat{V}$. 
Define a linear isomorphism $\hat{g}$ on $\hat{V}$ by
$\hat{g}(t^n\tensor a):= e^{-2\pii n} t^n\tensor ga$.
Then $\hat{g}$ defines an automorphism of a vertex algebra $\hat{V}$ 
and hence it gives rise to an automorphism of a Lie algebra
$\mathfrak{g}_V$.
Denote by $\mathfrak{g}_V^g$ the $\hat{g}$-invariants of
$\mathfrak{g}_V$, which is a Lie subalgebra of $\mathfrak{g}_V$.
Then the $g$-twisted universal enveloping algebra $\U^g(V)$ is defined 
to be the universal enveloping algebra for $\mathfrak{g}_V^g$.
The algebra $\U^g(V)$ has a universal property such that
for any weak $g$-twisted $V$-module $M$, the mapping  $a(n)\in
\U^g(V)\mapsto a_{(n)}= \res_z Y_M(a,z) z^n \in \End (M)$ gives 
a representation of $\U^g(V)$ on $M$.
It is clear that $\mathfrak{g}_V^g$ is spanned by images of elements 
$t^{n+\fr{r}{\abs{g}}}\tensor a$ with $a\in V^r$, $0\leq r\leq 
\abs{g}-1$.
We denote the image of $t^{n+\fr{r}{\abs{g}}}\tensor a$ in
$\mathfrak{g}_V^g$ by $a(n+\fr{r}{\abs{g}})$.
By definition, we have the following commutator relation:
\begin{equation}\label{commutator}
  [a(m), b(n)]=\sum_{i=0}^\infty \binom{m}{i} (a_{(i)}b)(m+n-i).
\end{equation}

\begin{df}
  For a monomial $x^1(n_1)\cds x^k(n_k)$ in $\U^g(V)$, we define its  
  {\it length} by $k$, {\it degree} by $\wt (x^1)+\cds +\wt (x^k)$ and 
  {\it weight} by $(\wt (x^1)-n_1-1)+\cds + (\wt (x^k)-n_k-1)$.
\end{df}

Let $W$ be a weak $g$-twisted $V$-module generated by one element
$w\in W$.  
In this case, a linear map $\phi_w : x^1(m_1)\cds x^k(m_k)\in
\U^g(V)\mapsto x^1_{(m_1)}\cds x^k_{(m_k)}w\in W=V\cd w$ gives a
surjection. 

The idea of the following assertion comes from M. Miyamoto 
\cite[Lemma 2.4]{M3}.

\begin{lem}\label{spanning set}
  Let $V$ be a $C_2$-cofinite VOA of CFT type and $W$ a weak
  $g$-twisted $V$-module generated by a non-zero element $w$,
  i.e. $W=V\cd w$.
  Let $U$ be a finite dimensional subspace of $V$ such that both
  $L(0)$ and $g$ act on $U$ and $V=U+C_2(V)$.
  Then the image $\phi_w(X)\in W$ of any monomial $X=x^1(m_1)\cds
  x^k(m_k)$ in $\U^g(V)$ can be expressed as a linear combination of
  images of monomials $\alpha^1(n_1)\cds \alpha^s(n_s)$ in $\U^g(V)$
  such that $\deg \alpha^1(n_1)\cds \alpha^s(n_s)$ is less or  
  equal to $\deg X$, $\wt \alpha^1(n_1)\cds \alpha^s(n_s)=\wt X$ and 
  $n_1<\cds <n_s<T$, where $T$ is a fixed element in
  $\fr{1}{\abs{g}}\Z$ such that $\phi_w(\beta (m))=0$ for all
  $\beta\in U$ and $m\geq T$.
\end{lem}

\pf
We divide the proof into several steps.

\begin{claim} 
  We can express the image $\phi_w (X)$ of any monomial $X=
  x^1(m_1)\cds x^k(m_k)\in \U^g(V)$ in the following form:
  $$
    \phi_w (X)= \phi_w(A) + \phi_w(B),
  $$
  where $A$ is a linear combination of monomials $\alpha^1(n_1)\cds
  \alpha^k(n_k)\in \U^g(V)$ with $\alpha^i \in U$ such that $\deg
  \alpha^1 (n_1)\cds \alpha^k(n_k) =\deg X$ and $\wt \alpha^1(n_1)\cds
  \alpha^k(n_k) =\alpha^k(n_k) =\wt X$, and $B$ is a sum of monomials
  whose degrees are less than $\deg X$ and weights are equal to 
  $\wt X$. 
\end{claim}

We prove the claim above by induction on $r=\deg X$.
The case $r=0$ is clear.
Assume that the claim is true for $r-1$.
Without loss, we may assume that both $L(0)$ and $g$ act on $x^i$, 
$1\leq i\leq k$ semisimply and none of them is the vacuum.
Then, by inductive assumption, $\phi_w(x^2(m_2)\cds x^k(m_k))$ can be
expressed a linear combination of images of monomials as stated.
Therefore, we may assume that $x^2,\cds,x^k$ are contained in $U$.
Since $V=U+C_2(V)$, we can write $x^1=\alpha^1+\sum_i a^i_{(-2)}b^i$
with $L(0)$-homogeneous $\alpha^1\in U$ and $a^i,b^i\in V$ such that
$\wt (\alpha^1)=\wt (a^i_{(-2)}b^i)=\wt (x^1)$.
Then $X=\alpha^1(m_1)x^2(m_2)\cds x^k(m_k) +\sum_i
(a^i_{(-2)}b^i)(m_1) x^2(m_2)\cds x^k(m_k)$.  
Then using \eqref{iterate} we can rewrite the image of second term in
the desired form because $\wt (a^i)+\wt (b^i)< \wt (a^i_{(-2)}b^i)$.
This completes the proof of Claim 1.

\begin{claim}
  Let $A=\alpha^1(m_1)\cds \alpha^k(m_k)\in \U^g(V)$ be a monomial
  with $\alpha^i\in U$ and $\sigma$ a permutation on the set $\{ 1, 2, 
  \dots, k\}$.
  Then we have the following equality in $W$:
  $$
    \phi_w\l( \alpha^{\sigma (1)}(n_{\sigma (1)})\cds \alpha^{\sigma
    (k)} (n_{\sigma (k)})\r) =\phi_w(A) + \phi_w(B),
  $$
  where $B$ is a sum of monomials whose degrees are less than $\deg A$
  and weights are equal to $\wt A$. 
\end{claim}

Again we proceed by induction on $r=\deg A$.
The case $r=0$ is obvious.
Assume that the assertion is correct for $\deg A= r-1$.
Then using the commutator formula \eqref{commutator} we can rearrange
$A$ to be as asserted since $\wt (\alpha^i_{(p)} \alpha^j)<\wt
(\alpha^i)+\wt (\alpha^j)$ for $p\geq 0$.
Thus, Claim 2 holds.

\begin{claim}
  Let $A= \alpha^1(m_1)\cds \alpha^k(m_k)\in \U^g(V)$ be a monomial
  with $\alpha^i\in U$ and $m_1\leq \cds \leq m_k<T$.
  Then the image $\phi_w(A)$ of $A$ can be expressed in the following
  form:  
  $$
    \phi_w(A) = \phi_w(B) +\phi_w(C), 
  $$
  where $B$ is a sum of monomials $\beta^1(n_1)\cds \beta^s(n_s)$ with   
  $\beta^j\in U$ such that $n_1<\cds <n_s$, $s\leq k$, $\deg
  \beta^1(n_1) \cds \beta^s(n_s)=\deg A$ and $\wt \beta^1(n_1)\cds
  \beta^s(n_s)=\wt A$, and $C$ is a sum of monomials whose degrees are 
  less than $\deg A$ and weights are equal to $\wt A$.
\end{claim}

We show that if the assertion is not correct then keeping both degree
and weight of $A$ we can make $m_1$ in a monomial $A$ infinitely
larger.  
We define an ordering on $\N\times \N$.
For $(r_1,s_1), (r_2,s_2)\in \N \times \N$, we define
$(r_1,s_1)>(r_2,s_2)$ if $r_1>r_2$, or $r_1=r_2$ and $s_1>s_2$.
By this ordering, $\N\times \N$ becomes a well-ordered set and hence
we can perform an induction on $(\deg A,\text{length} A)\in \N\times
\N$.
Clearly, the assertion is clear for $(\N,0)$, $(\N,1)$ and $(0,\N)$.
So we assume that the assertion is true for all elements in $\N\times
\N$ smaller than $(r,s)$ with $r>0$, $s>0$.
Then, by inductive assumption, we may assume that $m_2<\cds <m_k<T$. 
If $m_1<m_2$, then we are done.
So we have to consider the case $m_1=m_2$ and the case $m_1>m_2$.
But, the following argument shows that the latter case can be reduced
to the former case.
Assume that $m_1>m_2$. Then $A$ can be replaced by a linear
combination of $A'=\alpha^2(m_2)\alpha^1(m_1) \alpha^3(m_3)\cds
\alpha^k(m_k)$ and monomials whose degrees are smaller than $\deg A$ and
weights are the same as $\wt A$.
Then applying Claim 1 and Claim 2 together with inductive assumption
to $A'$, we can replace $A$ by a monomial $A''= (\alpha^1)'(m_1')\cds
(\alpha^k)'(m_k')$ such that $(\alpha^i)'\in U$, $m_1'>m_1$, $m_1'\leq
\cds \leq m_k'<T$, $\deg A''=\deg A$ and $\wt A''=\wt A$. 
Then, repeating this procedure, we will reach the case
$m_1=m_2<m_3<\cds <m_k$.

Now let us consider the case $m_1=m_2<m_3<\cds <m_k$.
In this case, both $\alpha^1$ and $\alpha^2$ are contained in  the
same eigenspace, say $V^r$. Write $m_1=n+\fr{r}{\abs{g}}$.
Using the iterate formula \eqref{iterate} on
$(\alpha^1_{(-1)}\alpha^2)_{(2n+1+\fr{2r}{\abs{g}})}$, we get
\begin{equation}\label{hoge}
\begin{array}{l}
  \ds
  \phi_w (\alpha^1(m_1)\alpha^2(m_1)\alpha^3(m_3)\cds \alpha^k(m_k))
  \vsb\\
  \qq = \lambda \phi_w((\alpha^1_{(-1)}\alpha^2)(2m_1+1)\alpha^3(m_3)\cds 
      \alpha^k(m_k)) 
  \vsb\\
  \qq  + \sum_{i>0} \mu_i \phi_w ( \alpha^1(m_1+i)\alpha^2(m_1-i)
      \alpha^3(m_3) \cds \alpha^k(m_k))
  \vsb\\
  \ds
  \qq
  + \sum_{i>0}\mu_i' \phi_w ( \alpha^2(m_1+i)\alpha^2(m_1-i)
    \alpha^3(m_3) \cds \alpha^k(m_k)) + \phi_w(X),
\end{array}
\end{equation}
where $X$ is a sum of monomials  whose degrees are less than $\deg A$
and weights are equal to $\wt A$.
Note that in the expansion of $(\alpha^1_{(-1)}\alpha^2)_{(2m_1
+1)}$, we can make the coefficient of $\alpha^1_{(m_1)}
\alpha^2_{(m_1)}$ non-zero by choosing suitable $N$ in
\eqref{iterate}.
The first term in the right-hand side of \eqref{hoge} has smaller
length than that of $A$ so that by induction together with Claim 1 and 
Claim 2 we may omit this term.
The second and third terms in the right-hand side of \eqref{hoge}
shall be reduced to the case $m_1>m_2$.
Therefore, we obtain a procedure which makes $m_1$ infinitely larger
with keeping $\deg A$ and $\wt A$, which must stop in finite steps.
Thus, we get Claim 3 and hence we complete the proof of the Lemma
\ref{spanning set}. 
\qed

\begin{rem}
  Even if $U$ is not finite dimensional, the lemma above still holds 
  when we can take a $T\in \fr{1}{\abs{g}}\Z$ such that $\phi_w(\beta(m))=0$ 
  for all $\beta\in U$, $m\geq T$.  
\end{rem}

\begin{rem}
  By Lemma \ref{spanning set} we can remove the repeat condition in
  \cite{B}.
\end{rem}

\begin{rem}
  There is another proof of Lemma \ref{spanning set} in \cite{NT}.
  See the proof of Theorem 3.2.7 of \cite{NT}.
\end{rem}

Now we can generalize Theorem 4.5 of \cite{ABD}.

\begin{thm}\label{g-regular}
  Every $g$-rational $C_2$-cofinite VOA of CFT type is $g$-regular.
\end{thm}

\pf
The proof is almost the same as that of Theorem 4.5 of \cite{ABD}.
The main idea in the proof of Theorem 4.5 of \cite{ABD} is to show
that every weak module has a non-trivial lowest weight vector.
By Lemma \ref{spanning set}, we can find a non-zero lowest weight
vector in every weak $g$-twisted module.
Thus, applying the argument in \cite{ABD} we get the assertion.
\qed

There are several corollaries of Lemma \ref{spanning set}.

\begin{cor}\label{cor1}
  Let $V$ be a $C_2$-cofinite VOA of CFT type.
  Then every irreducible weak $g$-twisted $V$-module $W$ is an
  irreducible ordinary $g$-twisted $V$-module.
\end{cor}

\pf
By Lemma \ref{spanning set}, we can introduce a
$\fr{1}{\abs{g}}\Z$-grading on $W$.
Therefore, every irreducible weak $g$-twisted module is exactly an
irreducible admissible $g$-twisted module.
Since every irreducible admissible module is an ordinary module,
we get the assertion.
\qed

\begin{cor}
  Let $V$ be a $C_2$-cofinite VOA of CFT type.
  Then every weak $g$-twisted $V$-module is admissible.
\end{cor}

\pf
By Proposition 3.6 of \cite{DLM2}, the $g$-twisted Zhu algebra
$A_g(V)$ (see \cite{DLM1}) is finite dimensional.
Then the argument in the proof of Proposition 5.6 of \cite{ABD} with
suitable modification leads to the assertion.
\qed

%\newpage

\section{Generalized theta functions on VOA-modules} %%%%%%%%%%%%

\subsection{Modular invariance of trace functions} %%%%%%%%%%%%

Let $V$ be a VOA and let $g$ and $h$ be mutually commutative
automorphisms on $V$. 
A $g$-twisted $V$-module $W$ is said to be {\it $h$-stable}
if there exists a linear isomorphism $\phi_W (h)$ on $W$ such that 
$$
  \phi_W (h) Y_W (a,z) = Y_W(ha,z)\phi_W (h)
$$
for all $a\in V$.
A linear isomorphism $\phi_W(h)$ is called {\it $h$-stabilizing
automorphism} or simply {\it stabilizing automorphism} on $W$.
For an ordinary $g$-twisted $h$-stable $V$-module $W$,
we can consider the following $q$-trace.
\begin{equation}\label{tr}
  T_W(a,\tau):=\tr_W z^{\wt (a)} Y_W(a,z) \phi_W(h)^{-1}
  q^{L(0)-c/24},
\end{equation}
where $q=e^{2\pii \tau}$ and $c$ denotes the central charge of $V$.
Zhu \cite{Z} proved that the space spanned by the trace functions
above is invariant under the action of the modular group $\SL_2(\Z)$
in the case $g=h=1$ and Dong, Li and Mason \cite{DLM2} generalized
his result to the case where $g$ and $h$ generate a finite abelian
subgroup in $\aut (V)$.  
Before we state their results, we have to introduce a structure
transformation of vertex operator algebras. 

\begin{df} (\cite[Theorem 4.2.1]{Z})
  Let $(V,Y(\cd,z),\vac,\w)$ be a VOA.
  For each homogeneous $a\in V$, the vertex operator 
  $$
    Y[a,z]:=e^{z\wt (a)} Y(a,e^z-1)=\sum_{n\in\Z} a_{[n]} z^{-n-1} 
    \in \End (V)[[z,z^{-1}]]
  $$
  provides another VOA structure on $V$ with the same vacuum vector
  $\vac$ and a new Virasoro vector $\tilde{\w}:=\w -(c/24)\vac$.
  We write $Y[\tilde{\w},z]=\sum_{n\in\Z} L[n]z^{-n-2}$ and denote
  $L[0]$-weight subspaces by $V_{[n]}= \{ a\in V\mid L[0]a=na\}$.
  Also, we use $\wt [a]$ to denote the $L[0]$-weight of $a\in V$.
\end{df}

We use the following notation.
For $\rho = \smat \in \SL_2(\Z)$, $(g,h)^\rho$ 
denotes $(g^\alpha h^\gamma, g^\beta h^\delta)$ and
$\rho \tau$ denotes $(\alpha \tau +\beta)(\gamma\tau +\delta)^{-1}$
for $\tau \in \{ z\in \C \mid \mathrm{Im}(z)>0\}$.

\begin{thm}\label{modular invariance}
  (\cite{DLM2})
  Let $V$ be a $C_2$-cofinite VOA and let $g$, $h$ be automorphisms on 
  $V$ generating a finite abelian subgroup in $\aut (V)$.
  Then the trace functions \eqref{tr} defined on irreducible
  $g$-twisted $h$-stable $V$-modules converge to linearly independent
  holomorphic functions on the  upper half plane. 
  Denote by $\mathcal{C}_1(g,h)$ the linear space spanned by the trace
  functions $T_W(a,\tau)$, where $W$ runs over irreducible
  $g$-twisted $h$-stable $V$-modules.
  For $\rho =\smat \in \SL_2(\Z)$ and $T_W(a,\tau)\in
  \mathcal{C}_1(g,h)$, define 
  \begin{equation}\label{action}
    (T_W)^\rho (a,\tau):= (\gamma \tau +\delta)^{-\wt [a]}
    T_W( a,\rho \tau ) .
  \end{equation}
  If $V$ is both $g$-rational and $g^\alpha h^\gamma$-rational, then  
  $\rho$ defines a linear isomorphism from $\mathcal{C}_1(g,h)$
  to $\mathcal{C}_1((g,h)^\rho)$.
\end{thm}

The space $\mathcal{C}_1(g,h)$ is called a {\it genus one twisted 
conformal block}.

\begin{rem}
  This theorem has been generalized to involve intertwining operators
  in \cite{M2} and \cite{Y}.
\end{rem}

\subsection{Semisimple primary vectors} %%%%%%%%%%%%%%%%%%%%%%%%%

Here we review the theory of the physics superselection principal and
semisimple primary vectors in \cite{L1}.

\begin{df}
  A vector $u\in V$ is called {\it a semisimple primary vector} if it
  satisfies the following.
  \vsb\\
  (i)\ $L(n)u=\delta_{n,0}h$ for $n\geq 0$.
  \vsb\\
  (ii)\ $u_{(m)}u =\delta_{m,1}\gamma \vac$ for $m\geq 0$ and some   
    $\gamma\in\Q$. 
  \vsb\\
  (iii)\ $u_{(0)}$ acts on $V$ semisimply.
\end{df}

Since $u_{(0)}$ is a derivative operator and keeps each homogeneous
subspace of $V$, its exponential operator $\exp (\alpha u_{(0)})$
gives an automorphism of $V$ for any $\alpha \in\C$. 
In the following, we denote $\exp(-2\pii u_{(0)})$ by $\sigma (u)$.
If all eigenvalues of $u_{(0)}$ on $V$ is contained in $\fr{1}{T}\Z$
for some $T\in\Z$, then $\sigma (u)$ have a finite order.
We call such a semisimple primary vector {\it rational}.

Let $u$ be a rational semisimple primary vector and $g$ an
automorphism of $V$ of finite order such that $gu=u$ (so $g\sigma (u) 
=\sigma (u) g$).
%We will freely use $u(n)$ to denote $u_{(n)}$ for semisimple primary
%vector $u$.
Define 
\begin{equation}\label{delta-operator}
  \Delta (u,z)
  := z^{u_{(0)}} \exp \l( -\sum_{n=1}^\infty \fr{u_{(n)}}{n}
     (-z)^{-n}\r) .
\end{equation}
Since $u_{(0)}$ acts on $V$ semisimply, $\D (u,z)$ is a well-defined  
operator on $V$.
Let $(W,Y_W(\cd,z))$ be a weak $g$-twisted $V$-module.
The following proposition is due to Li \cite{L2} (see also
\cite{L1}). 

\begin{prop}\label{transform}
  (\cite[Proposition 5.4]{L2})\ 
  $(W,Y_W(\D (u,z)\,\cd,z))$ is a weak $g\sigma (u)$-twisted
  $V$-module. 
\end{prop}

Let us denote $(W,Y_W(\D (u,z)\,\cd,z))$ simply by $\tilde{W}$.
We can write the action of $a\in V$ on $\tilde{W}$ in the
following way.
By the assertion above, there exists a linear isomorphism $\varphi_W : 
W \to \tilde{W}$ such that $Y_{\tilde{W}}(a,z)\varphi_W = \varphi_W
Y_W(\D (u,z)a,z)$ for all $a\in V$.
Define the Schur polynomials $p_s(x_1,x_2,\dots)$ in variables 
$x_1, x_2,\dots$ by the following equation:
\begin{equation}\label{Schur polynomial}
  \exp \l( -\dsum_{n=1}^\infty \dfr{x_n}{n} (-z)^{-n}\r)
  = \dsum_{s=0}^\infty p_s(x_1,x_2,\dots) z^{-s}.
\end{equation}
Assume that $u_{(0)}a=\lambda a$ for some $\lambda\in\Q$.
%For $a\in V$, write $Y_{\tilde{W}}(a,z)=\sum_{n\in\Q}
%\tilde{a}_{(n)}z^{-n-1}$. 
Then the vertex operator of $a$ on $\tilde{W}$ is given as
follows:
\begin{equation}\label{concrete}
  Y_{\tilde{W}}(a,z)\varphi_W =\varphi_W Y_W(\D (u,z)a,z)
  = \varphi_W \sum_{s=0}^\infty z^{-s+\lambda} 
    Y_W(p_s(u_{(1)},u_{(2)},\dots) a,z).
\end{equation}

%In particular, we have $L(0)\varphi_W =\varphi_W
%(L(0)+u(0)+\fr{1}{2}\gamma)$. 

The Delta operator $\D (u,z)$ has an additive property.
A pair of semisimple primary vectors $u$ and $v$ such that
$u_{(0)}v=0$ is called {\it mutually commutative} because
we have  $\D (u,z) \D (v,z)=\D (v,z) \D (u,z)= \D (u+v,z)$.
In particular, $u$ is commutative with itself so that $\D (u,z)$ is
invertible because $\D (u,z) \D (-u,z)= \D (0,z)=\id_V$.
The following statement is easy.

\begin{prop}\label{correspondence}
  We have a bijective correspondence between the set of 
  irreducible $g$-twisted $V$-modules and the set of irreducible
  $g\sigma (u)$-twisted modules through the Delta operator $\D (u,z)$.
  Furthermore, if an automorphism $h$ on $V$ is commutative with $g$
  and acts on $u$ trivially, then the set of irreducible $g$-twisted
  $h$-stable $V$-modules and the set of irreducible $g\sigma
  (u)$-twisted $h$-stable $V$-modules are in one-to-one
  correspondence. 
\end{prop}

We will need the following lemma.

\begin{lem}\label{translation}
  Let $V$ be a $g$-rational $C_2$-cofinite VOA of CFT type.
  Then $V$ is $g\sigma (u)$-regular for every rational semisimple
  primary vector $u$.
\end{lem}

\pf
In this case, $V$ is $g$-regular by Theorem \ref{g-regular}.
Let $(W,Y_W(\cd,z))$ be a weak $g\sigma (u)$-twisted $V$-module.
Then $(W,Y_W(\D (-u,z)\,\cd,z))$ is a weak $g$-twisted $V$-module.
Since $V$ is $g$-regular,  $(W,Y_W(\D (-u,z)\,\cd,z))$ is a direct sum 
of irreducible $g$-twisted $V$-modules. 
Then $(W,Y_W(\D (u,z) \D(-u,z)\,\cd,z))=(W,Y_W(\cd,z))$ is also a 
direct sum of irreducible $g\sigma (u)$-twisted $V$-modules.
\qed

\subsection{Main Theorem} %%%%%%%%%%%%%%%%%%%%%%%%%%%%%%%%%%%%%%%

In the following context, we will work over the following setting.
\vsb\\
(1)\ $V$ is a $C_2$-cofinite vertex operator algebra of CFT type.
\vsb\\
(2)\ $L(1)V_1=0$. 
\vsb\\
(3)\ $g$ and $h$ are automorphisms on $V$ generating a finite
abelian subgroup in $\aut (V)$.
\vsb\\
(4)\ $V$ is $k$-rational for all $k\in \la g,h\ra$.
\vsb\\
(5)\ $H$ is a set of mutually commutative rational semisimple primary
vectors in $V^{\la g,h\ra}$, where $V^{\la g,h\ra}$ denotes the fixed
point sub VOA under $\la g,h\ra$.  
\vsb\\
(6)\ For $u\in H$, $\sigma (u)$ denotes $\exp (-2\pii u_{(0)})$ and
$\lambda_u(a)$ is a linear function on $V$ defined as
$u_{(0)}a=\lambda_u(a)a$ for $a\in V$.
\vsb\\
(7)\ $E(H)=\{ \sigma (u)\mid u\in H\}$, an abelian subgroup of $\aut 
(V)$. 
\vsb\\
We make some remarks on the assumption above.
By (2), $V$ possesses the unique invariant bilinear form $\la \cd
,\cd \ra$ such that $\la \vac,\vac\ra =-1$ (cf. \cite{L3}).
Note that this bilinear form satisfies $a_{(1)}b=\la a,b\ra \vac$ for
$a,b\in V_1$. 
The assumption (4) is satisfied if $V^{\la k\ra}$ is rational for all
$k\in \la g,h\ra$ by Proposition \ref{g-rationality}.
All $\sigma (ru)$ with $u\in H$ are finite automorphisms on $V$ for
any $r\in \Q$.  Therefore, $H$ forms a $\Q$-vector space in $V_1$.
Since $H$ is contained in $V^{\la g,h\ra}$, we have $[\la g,h\ra,
E(H)]=1$ in $\aut (V)$. 

Let $k_1,k_2\in \la g,h\ra$.
Since $V$ is $k_1$-rational, there are finitely many irreducible
inequivalent $k_1$-twisted $k_2$-stable $V$-modules.
We denote the complete set of inequivalent irreducible $k_1$-twisted 
$k_2$-stable $V$-modules by $\{ (W^i(g,h),\phi_{i,k_1}(k_2)) \mid 
i=1,2,\dots,N=N(k_1,k_2)\}$, where $\phi_{i,k_1}(k_2)$ are
(fixed) $k_2$-stabilizing automorphisms on $W^i(k_1,k_2)$.
Note that the number $N(k_1,k_2)$ of irreducible $k_1$-twisted $k_2$-stable
$V$-modules is the same as that of irreducible $k_1^\alpha
k_2^\gamma$-twisted $k_1^\beta k_2^\delta$-stable $V$-modules for
all $\binom{\alpha\ \beta}{\gamma\ \delta}\in \SL_2(\Z)$ 
by Theorem \ref{modular invariance}.
Recall the genus one twisted conformal block $\mathcal{C}_1(g,h)$ 
which is a linear span of $q$-traces
$$
  T_{W^i(g,h)}(a,\tau)=\tr_{W^i(g,h)} Y(a,z) 
    \phi_{i,g}(h)^{-1} q^{L(0)-c/24},\ 1\leq i\leq N.
$$
For $\rho =\smat \in \SL_2(\Z)$, we have the following
transformation:
\begin{equation}\label{transformation}
  T_{W^i(g,h)}(a,\rho \tau)
  = (\gamma \tau +\delta)^{\wt [a]} 
    \sum_{j=1}^{N} A_{ij}(\rho,(g,h)) 
    T_{W^j((g,h)^\rho)} (a,\tau),
\end{equation}
where the constants $A_{ij}(\rho,(g,h))$ are given by Theorem
\ref{modular invariance} and independent of $a$ and $\tau$. 

\begin{df}
  For $u,v\in H$ and $a\in V$, define
  \begin{equation}
    \begin{array}{l}
      Z_{W^i(g,h)}(a;(u,v);\tau)
      \vsb\\
      \qq := \dsum_{s=0}^\infty \tr_{W^i(g,h)}
        \Big\{ \l( p_s(u_{(1)},u_{(2)},\dots)a\r)_{(\wt (a) +
        \lambda_u (a) -s-1)}
      \vsb\\
      \qq\qq  \times e^{\pii \la u,v\ra}
        \exp \l( 2\pii v_{(0)} \r) \phi_{i,g}(h)^{-1} 
        q^{L(0)+u_{(0)} +\fr{1}{2}\la u,u\ra -c/24}\Big\} ,
    \end{array}
  \end{equation}
  where $q$ denotes $e^{2\pii\tau}$, $p_s(x_1,x_2,\dots)$ is the Schur
  polynomial defined by \eqref{Schur polynomial} and $\lambda_u(a)$ is
  a scalar such that $u_{(0)}a=\lambda_u(a)u$. 
  We call $Z_{W^i(g,h)}(a;(u,v);\tau)$ {\it a generalized theta
  function on $W^i(g,h)$ with respect to $H$}.
\end{df}

\begin{rem}
  In \cite{M1}, Miyamoto defined the function above in the case when
  $g=h=1$ and $a=\vac$ and he called
  $Z_{W^i(1,1)}(\vac;(0,v);\tau) \eta (\tau)^c$ {\it a theta function
  of $W^i(1,1)$}, where $\eta (\tau)=q^{1/24} \prod_{n=1}^\infty
  (1-q^n)$ is the Dedekind eta function.
\end{rem}

We consider the modular transformations of $Z_{W^i(g,h)}(a;(u,v);\tau)$. 
Our main theorem is the following.

\begin{thm}\label{main thm}
  The generalized theta function $Z_{W^i(g,h)}(a;(u,v);\tau)$
  converges to a holomorphic function on the upper half plane and 
  gives a vector in $\mathcal{C}_1(g\sigma (u),h\sigma (v))$ 
  for each $1\leq i\leq N$.
  Furthermore, we have the following modular transformation for
  $\rho=\smat \in \SL_2(\Z) :$
  $$ %\begin{equation} \label{main}
  \begin{array}{l}
    \ds
    (\gamma \tau +\delta )^{-\wt [a]} Z_{W^i(g,h)}
     \l( a;(u,v);\dfr{\alpha\tau +\beta}{\gamma \tau +\delta} \r)
    \vsb\\
    \ds
    \q = 
    %(\gamma\tau +\delta )^{\wt [a]} 
    \sum_{i=1}^N A_{ij}(\rho,(g,h)) \
    Z_{W^i((g,h)^\rho)}(a;(\alpha u+\gamma v,\beta u+\delta
    v);\tau) ,
  \end{array}
  $$ %\end{equation}
  where $A_{ij}(\rho,(g,h))$ are the constants given by \cite{DLM2} 
  as in the equation \eqref{transformation}.
\end{thm}

\begin{rem}
  Since $L[0]=L(0)+\sum_{i\geq 1}c_i L(i)$ for some $c_i\in \C$, 
  we have $[L[0],u_{(0)}]=0$ for all $u\in H$.
  Therefore, $u_{(0)}$ acts on each $L[0]$-homogeneous subspace
  $V_{[n]}$ diagonally.
\end{rem}

\pf
We divide the proof into two parts.
In the first part, we show that the modular transformation 
$(\gamma \tau +\delta)^{-\wt [a]} Z_{W^i(g,h)}(a;(u,v);\rho \tau)$ 
is uniquely expressed as a linear combination of
$Z_{W^j((g,h)^\rho)}(a;(\alpha u+\gamma v,\beta u+\delta
v);\tau)$, $1\leq j\leq N$. 
Then we show that the coefficients of the linear combination are
exactly given as stated.

Let us consider the meaning of $Z_{W^i(g,h)}(a;(u,v);\rho \tau)$.
By Proposition \ref{transform}, $\tilde{W}^i(g,h)$ $:=
(W^i(g,h),Y(\D (u,z)\,\cd,z))$ is an irreducible $g \sigma
(u)$-twisted $V$-module. 
Since $H$ is a subspace of $V^{\la g,h\ra}$, $v_{(0)}$ acts on
$\tilde{W}^i(g,h)$ as a derivation
$$
  v_{(0)} Y_{\tilde{W}^i(g,h)} (a,z)
  = Y_{\tilde{W}^i(g,h)} (v_{(0)}a,z) 
    + Y_{\tilde{W}^i (g,h)} (a,z) v_{(0)}.
$$
Therefore, we have the following.
$$
  \exp \l( -2\pii v_{(0)}\r) Y_{\tilde{W}^i(g,h)}(a,z)
  = Y_{\tilde{W}^i(g,h)} (\sigma (v)a,z) 
    \exp \l( -2\pii v_{(0)}\r) .
$$
(Note that $\exp \l( -2\pii v_{(0)}\r)$ is well-defined on
$\tilde{W}^i(g,h)$ since it is irreducible.)
Namely, $\tilde{\sigma}(v):= \exp \l( -2\pii v_{(0)}\r)$, here
$v_{(0)}=\res_z Y_{\tilde{W}^i(g,h)}(v,z) \in \End (\tilde{W}^i(g,h))$, 
is a $\sigma (v)$-stabilizing automorphism on $\tilde{W^i}(g,h)$.
On the other hand, by definition, there exist linear isomorphisms 
$\varphi_{W^i(g,h)} : W^i(g,h)\to \tilde{W}^i(g,h)$, 
$1\leq i\leq N$, such that 
$$
  Y_{\tilde{W}^i(g,h)}(a,z)\varphi_{W^i(g,h)} 
  = \varphi_{W^i(g,h)} Y_{W^i(g,h)}(\D (u,z)a,z) .
$$
Then 
$$
\begin{array}{l}
  \varphi_{W^i(g,h)} \phi_{i,g}(h)
    \varphi_{W^i(g,h)}^{-1} Y_{\tilde{W}^i(g,h)}(a,z) 
  \vsb\\
  = \varphi_{W^i(g,h)} \phi_{i,g}(h) 
    Y_{W^i(g,h)}(\D (u,z)a,z) \varphi_{W^i(g,h)}^{-1}
  \vsb\\
  = \varphi_{W^k(g,h)} Y_{W^i(g,h)}(\D (u,z)h a,z) 
    \phi_{i,g}(h)\varphi_{W^i(g,h)}^{-1}
  \vsb\\
  = Y_{\tilde{W}^i(g,h)}(h a,z) \varphi_{W^i(g,h)}
  \phi_{i,g}(h) \varphi_{W^i(g,h)}^{-1}.
\end{array}
$$
Thus, a composition $\tilde{\phi}_{i,g}(h):= \varphi_{W^i(g,h)}
\phi_{i,g}(h) \varphi_{W^i(g,h)}^{-1}$ provides a 
$k_2$-stabilizing automorphism on $\tilde{W}^i(g,h)$ for each
$1\leq i\leq N$.  
Therefore, we see that all inequivalent irreducible $g \sigma (u)$-twisted 
$h \sigma (v)$-stable $V$-modules are filled by 
$\tilde{W}^i(g,h)$, $1\leq i\leq N$ with the stabilizing
automorphisms $\tilde{\phi}_{i,g}(h)\tilde{\sigma} (v)$ by
Proposition \ref{correspondence}. 
Hence, by Theorem \ref{modular invariance}, the trace function
$$
  T_{\tilde{W}^i(g,h)}(a;(u,v);\tau ):=
  \tr_{\tilde{W^i}(g,h)} z^{\wt (a)} Y_{\tilde{W}^i(g,h)}
  (a,z) \tilde{\sigma}(v)^{-1} \tilde{\phi}_{i,g} (h)^{-1}
  q^{L(0)-c/24} 
$$
converges on the upper half plane and gives a vector in
$\mathcal{C}_1(g \sigma (u),h \sigma(v))$.
Since $V$ is $g \sigma (u)$-rational by Lemma \ref{translation},
$\mathcal{C}_1(g \sigma (u),h \sigma (v))$ is spanned by
$T_{\tilde{W}^i(g,h)} (a;(u,v);\tau)$, $1\leq i\leq N$.
Therefore, by Theorem \ref{modular invariance},
$\rho =\smat \in \SL_2(\Z)$ defines a linear isomorphism between 
$\mathcal{C}_1(g \sigma (u),h \sigma (v))$ and 
$\mathcal{C}_1((g \sigma(u), h \sigma (v))^\rho)$ in the following
way: 
\begin{equation}\label{hoge2}
\begin{array}{l}
  \ds T_{\tilde{W}^i(g,h)} (a;(u,v);\rho\tau )
  \vsb\\ \ds
  \q = (\gamma \tau +\delta)^{\wt [a]} 
    \sum_{j=1}^N B_{ij}((g,h),(u,v),\rho) 
    T_{\tilde{W}^i((g,h)^\rho)} (a;(\alpha u+\gamma v,\beta
    u+\delta v);\tau),
\end{array}
\end{equation}
where $B_{ij}((g,h),(u,v),\rho)$ are scalars independent of $a$ 
and $\tau$. 

By the way, using $Y_{\tilde{W}^i(g,h)}(a,z)
\varphi_{W^i(g,h)} = \varphi_{W^i(g,h)} Y_{W^i(g,h)}(\D (u,z)a,z)$, 
we have
$$
\begin{array}{l}
  T_{\tilde{W}^i(g,h)}(a;(u,v);\tau )
  \vsb\\
  = \tr_{\tilde{W}^i(g,h)} z^{\wt (a)} Y(a,z) 
    \tilde{\sigma} (v)^{-1} \tilde{\phi}_{i,g}(h)^{-1}
    q^{L(0)-c/24} 
  \vsb\\
  = \tr_{W^i(g,h)} z^{\wt (a)} Y(a,z) 
    \tilde{\sigma} (v)^{-1} \tilde{\phi}_{i,g}(h)^{-1}
    q^{L(0)-c/24} \varphi_{W^i(g,h)}
  \vsb\\
  = \tr_{W^i(g,h)} \varphi_{W^i(g,h)} z^{\wt (a)} 
    Y(\D (u,z)a,z) \exp\l( 2\pii (v_{(0)}+\la u,v\ra )\r)
  \vsb\\
  \hspace{2cm} \times
    \phi_{i,g}(h)^{-1} q^{L(0)+u_{(0)}+\fr{1}{2}\la u,u\ra -c/24} 
  \vsb\\
  = e^{\pii \la u,v\ra} Z_{W^i(g,h)} (a;(u,v);\tau),
\end{array}
$$
where we have also used that
$v_{(0)} \varphi_{W^i(g,h)} = \varphi_{W^i(g,h)} (v_{(0)}+\la 
u,v\ra)$ and $L(0) \varphi_{W^i(g,h)} = \varphi_{W^i(g,h)}
(L(0)+u_{(0)}+\fr{1}{2}\la u,u\ra )$. 
Therefore, by \eqref{hoge2}, $(\gamma \tau +\delta)^{-\wt [a]}
Z_{W^i(g,h)}(a;(u,v);\rho \tau)$ is a linear combination of 
$Z_{W^j((g,h)^\rho)}(a;(\alpha u+\gamma v,\beta u+\delta
v);\tau)$.

Next, we show that $B_{ij}((g,h),(u,v),\rho) = 
e^{\pii (-\la \alpha u+\gamma v,\beta u+\delta v\ra +\la u,v\ra)} 
A_{ij} (\rho,(g,h))$ for $1\leq i,j\leq N$, which would complete the proof. 
Since $B_{ij}((g,h),(u,v),\rho)$, $1\leq i,j\leq N$, are independent 
of $a$, we prove the equality in the case where $a=\vac$.
We use some results from \cite{Z} and \cite{M1}.

For $a_1,\dots,a_n\in V^{\la g,h\ra}$, set 
$$
\begin{array}{l}
  S_{W^i(g,h)} ((a_1,z_1),\dots,(a_n,z_n);\tau)
  \vsb\\
  := q_{z_1}^{\wt (a_1)} \cds q_{z_n}^{\wt (a_n)} 
     \tr_{W^i(g,h)} Y(a_1,q_{z_1})\cds Y(a_n,q_{z_n})
     \phi_{i,g}(h)^{-1}q^{L(0)-c/24}_{\tau},
\end{array}
$$
where $q_x$ denotes $e^{2\pii x}$.
We deduce a recurrent formula for $S_{W^i(g,h)}$.
Before we state it, we introduce the following functions (cf.\ 
\cite{Z} and \cite{DLM2}).

The Eisenstein series $G_{2k}(\tau)$ $(k=2,3,\dots)$ are series
\begin{equation}\label{Eisenstein}
  G_{2k}(\tau) = 2\zeta (2k) +\dfr{2(2\pii)^{2k}}{(2k-1)!}
  \dsum_{n=1}^\infty \dfr{n^{2k-1}q^n}{1-q^n},
\end{equation}
where $\zeta (2k)=\sum_{n=1}^\infty 1/n^{2k}$.
We use normalized Eisenstein series
\begin{equation}
  E_{2k}(\tau):= \dfr{1}{(2\pii)^{2k}} G_{2k}(\tau).
\end{equation}
Since $G_{2k}(\tau)$ is a modular form of weight $2k$ for the modular
group $\SL_2(\Z)$, we have
\begin{equation}\label{moge1}
  E_{2k}\l(\dfr{\alpha\tau +\beta}{\gamma \tau +\delta}\r)
  = (\gamma \tau +\delta)^{2k} E_{2k}(\tau)\q \text{for}\ 
  \mat\in \SL_2(\Z).
\end{equation}

We define the functions $\wp_k(z,\tau)$ $(k\geq 1)$ by
\begin{equation}
  \wp_k(z,\tau)
  := \dfr{1}{z^k}+(-1)^k \dsum_{n=1}^\infty \binom{2n+1}{k-1}
     G_{2n+2}(\tau)z^{2n+2-k}.
\end{equation}
The following modular transformations are well-known:
\begin{equation}\label{moge2}
  \wp_k \l(\dfr{z}{\gamma \tau +\delta},
   \dfr{\alpha \tau +\beta}{\gamma \tau +\delta}\r)
   =(\gamma \tau +\delta)^k \wp_k(z,\tau)\q \text{for}\ 
   \mat\in \SL_2(\Z).
\end{equation}

\begin{ass}\label{recurrent}
  (\cite[Proposition 4.4.2]{Z})\ 
  For $a_1,\dots,a_n, b\in V^{\la g,h\ra}$, the following recurrent
  formula holds: 
  $$
  \begin{array}{l}
     S_{W^i(g,h)} 
       \l( (b,x),(a_1,z_1),\dots,(a_n,z_n);\tau \r)  
     \vsb\\
     = S_{W^i(g,h)} 
       \l( (b_{[-1]}a_1,z_1),(a_2,z_2),\dots,(a_n,z_n);\tau\r) 
     \vsb\\
     \qq 
     - \dsum_{k=2}^\infty E_{2k}(\tau ) S_{W^i(g,h)} 
       \l( (b_{[2k-1]}a_1,z_1), (a_2,z_2),\dots,(a_n,z_n);\tau\r)
     \vsb\\
     \qq 
     + \dsum_{s=1}^n \dsum_{m=0}^\infty \dfr{1}{(2\pii )^{m+1}}
       \l\{ \wp_{m+1}(x-z_s,\tau)-\wp_{m+1}(z_1-z_s,\tau)\r\}
     \vsb\\
     \qq\qq \times
      S_{W^i(g,h)} 
      \l( (a_1,z_1),\dots,(b_{[m]}a_s,z_s),\dots,(a_n,z_n);\tau\r) .
  \end{array}
  $$
\end{ass}

\pf
Since $W^i(g,h)$ are untwisted ordinary modules for $V^{\la
g,h\ra}$,  we can use the same argument as that in \cite{Z} and hence
we obtain the same consequence.
(Note that our usage of the Eisenstein series differs from that of
Zhu in \cite{Z} by scalar multiples.)
\qed

Using the recurrent formula above, we can show the following.

\begin{ass}
  For $a_1,\dots,a_n\in V^{\la g,h\ra}$ and $\rho =\smat \in
  \SL_2(\Z)$, we have  
  $$
  \begin{array}{l}
  \ds
  S_{W^i(g,h)}\l( \l( a_1,\fr{z_1}{\gamma \tau+\delta}\r) ,\dots, 
    \l( a_n,\fr{z_n}{\gamma\tau +\delta}\r); \fr{\alpha\tau
    +\beta}{\gamma \tau +\delta} \r)
  \vsb\\
  \q = (\gamma\tau+\delta)^{\wt [a_1]+\cds +\wt [a_n]} 
    \dsum_{j=1}^N A_{ij}(\rho,(g,h)) S_{W^j((g,h)^\rho)} 
    ((a_1,z_1),\dots,(a_n,z_n);\tau),
  \end{array}
  $$ 
  where scalars $A_{ij}(\rho,(g,h))$ are given by
  \eqref{transformation}. 
\end{ass}
  
\pf
We proceed by induction on $n$.
The case $n=1$ is already known by Theorem \ref{modular invariance}. 
Since the $n$-point trace $S_{W^i(g,h)}$ is completely determined
by the $1$-point trace $T_{W^i(g,h)}$, using \eqref{moge1},
\eqref{moge2} and the recurrent formula obtained in Assertion
\ref{recurrent}, we get the assertion. 
\qed

Now, by the assertion above, we are in position to use the same
argument as that in \cite{M1}. 
We show the following brief generalization.

\begin{ass}\label{a=vac}
  (\cite[Main Theorem]{M1})
  \\
  Set $Z_{W^i(g,h)}'(u;v;\tau):=
    Z_{W^i(g,h)}(\vac;(u,v);\tau)$.  
  For $\rho =\smat \in \SL_2(\Z)$, we have
  \begin{equation}\label{moge3}
  \begin{array}{l}
    \ds
    Z'_{W^i(g,h)} \l( u;v; \fr{\alpha\tau +\beta}{\gamma \tau
      +\delta} \r)
    \vsb\\
    \ds
    \qq = \sum_{j=1}^N A_{ij}(\rho,(g,h)) \
      Z'_{W^j((g,h)^\rho)} 
      (\alpha u + \gamma v; \beta u+\delta v;\tau ).
  \end{array}
  \end{equation}
\end{ass}

\pf
The formula above was proved when $g=h=1$ in \cite{M1}.
The proof in \cite{M1} is given by direct calculations on the
$n$-point traces.
So by tracing calculations in \cite{M1} together with Theorem
\ref{modular invariance} and Assertion 2, one can verify the assertion.
\qed

Therefore, we have reached $B((g,h),(u,v),\rho)= 
e^{\pii (-\la \alpha u+\gamma v,\beta u -\delta v\ra +\la u,v\ra)} 
A((g,h),\rho)$ by \eqref{moge3}.
This completes the proof of Theorem \ref{main thm}.
\qed

There is an interesting consequence of Theorem \ref{main thm}.
For each $\rho\in \SL_2(\Z)$, let us denote by
$\Psi_{(g,h)}(\rho)$ the isomorphism from $\mathcal{C}_1(g,h)$
to $\mathcal{C}_1((g,h)^\rho)$ given as \eqref{action}.
The proof of Theorem \ref{main thm} tells us that the space
$\mathcal{C}_1(g \sigma (u),h \sigma (v))$ is spanned by functions
$Z_{W^i(g,h)}(a;(u,v);\tau)$, $1\leq i\leq N$ and the matrix
representation of $\Psi_{(g,h)}(\rho)$ and that of
$\Psi_{(g \sigma (u),h \sigma (v))} (\rho)$ are given by the same 
matrix $(A_{ij}(\rho,(g,h)))_{ij}$. 
Namely, we have proved that the internal automorphisms do not 
change the genus one twisted conformal blocks:

\begin{cor}\label{main cor}
  For a pair $(u,v)\in H\times H$, define a linear isomorphism
  $\Omega_{(g,h)}(u,v):\mathcal{C}_1(g,h) \to
  \mathcal{C}_1(g \sigma (u),h \sigma (v))$ 
  by $T_{W^i(g,h)}(a ,\tau) \mapsto Z_{W^i(g,h)}
  (a;(u,v);\tau)$.   
  Then we have 
  $$
    \Psi_{(g \sigma (u), h \sigma (v))}(\rho) \circ 
    \Omega_{(g,h)}(u,v)
    = \Omega_{(g,h)^\rho}(\alpha u +\gamma v,\beta u+\delta v)
      \circ \Psi_{(g,h)}(\rho) 
  $$
  for every $\rho = \smat \in \SL_2(\Z)$.
\end{cor}

%We also note that a subgroup $E(H)$ of $\aut (V)$ has a special
%importance in \cite{DG}. 

\subsection{Relation to abelian coset construction}

We keep the setup of the previous subsection.
In this subsection we consider the case where $g=h=1$ and $V$ is simple.
In addition, we assume that the restriction of the invariant bilinear form 
$\la \cd , \cd \ra$ on $H$ is non-degenerate.
Set $\h=\C \tensor_\Q H \subset V_1$ and let $\{ h^1,\dots,h^{\dim \h}\}$
be an orthonormal basis of $\h$.
For $h,k \in \h$, their vertex operators satisfies the commutator
relation $[h_{(m)},k_{(n)}]= \delta_{m+n,0} m \la h,k\ra$. 
Therefore, $\h$ generates a free bosonic sub VOA $M_\h (1,0)$ with the 
Virasoro vector 
$$
  \w_\h := \hf \dsum_{i=1}^{\dim \h} h_{(-1)}^ih_{(-1)}^i \vac .
$$
For $\alpha \in \h$, set 
$$
  V^\alpha :=\{ x\in V \mid h_{(0)}x=\la \alpha,h \ra x\q \mbox{for }
  \ h\in \h\} ,
$$
and we define $L=\la \alpha \in \h \mid V^\alpha \ne 0\ra$ which 
is a subgroup of the additive group $\h$.
Then we obtain an $L$-graded structure 
$$
  V=\oplus_{\alpha \in L} V^\alpha 
$$
with $Y(a,z) b \subset V^{\alpha+\beta}((z))$ for $a\in V^\alpha$, 
$b\in V^\beta$.
Note that $V^0$ is a sub VOA.
Since we have assumed that $V$ is $C_2$-cofinite, $L$ is a finitely
generated free $\Z$-module. Namely, $L$ equipped with $\la\cd,\cd\ra$ 
is a rational lattice.
As $V$ is simple,  $L=\{ \alpha \in \h \mid V^\alpha \ne 0\}$.
Define the space of highest weight vectors 
$$
  \Omega_V:= \{ x\in V \mid h_{(n)} x=0 \q \mbox{for }\ h\in \h,\ 
  n\geq 1\} .
$$
Since $[h_{(0)}, k_{(m)}]=0$ for $h,k\in \h$, $m\in \Z$, $h_{(0)}$ 
preserves $\Omega_V$.
Then 
$$
  \Omega_V =\oplus_{\alpha \in L} \Omega_V^\alpha,\q \mbox{where }\ 
  \Omega_V^\alpha = \Omega_V\cap V^\alpha.
$$
Since $V\simeq M_\h(1,0)\tensor_\C \Omega_V$ as a linear space, we 
have the following decomposition:
$$
  V=M_\h(1,0)\tensor \Omega_V 
  = \oplus_{\alpha \in L} M_\h(1,0)\tensor \Omega_V^\alpha.
$$
Note that $\Omega_V^0$ is a sub VOA of $V^0$ with the Virasoro vector 
$\w_\Omega:=\w -\w_\h$ and is a commutant subalgebra of $M_\h(1,0)$ 
such that $V^0=M_\h(1,0)\tensor \Omega_V^0$.
For $a\in \Omega_V^\alpha$, define
$$
  Y_{\Omega_V}(a,z):= E^-(\alpha,z) Y_V(a,z) E^+(\alpha,z) z^{-\alpha_{(0)}},
$$
where
$$
  E^\pm(\alpha,z):=\exp \l( \dsum_{n=1}^\infty \dfr{\alpha_{(\pm n)}}{\pm n}
  z^{\mp n}\r) .
$$
Then it is shown in \cite{DL} and \cite{L4} that the structure $(\Omega_V,
Y_{\Omega_V}(\cd,z),\vac,\w_\Omega)$ is a simple $L$-graded generalized 
vertex operator algebra  with central charge $c-\dim \h$.

Let $W$ be an irreducible $V$-module.
Then we also have a decomposition
$$
  W=M_\h(1,0)\tensor \Omega_W ,
$$ 
where $\Omega_W:= \{ w\in W \mid h_{(n)} w=0\ \mbox{for}\ h\in \h,\ 
n\geq 1\}$ denotes the space of highest weight vectors in $W$.
For $\lambda \in \h$, set $W^\lambda = \{ w\in W \mid h_{(0)} w
=\la \lambda, h\ra w\ \mbox{for}\ h\in \h \}$. 
Since $W$ is irreducible, there is an $L$-subset $\lambda_W +L$ of 
$\h$ such that $W=\oplus_{\beta\in L+\lambda_W} W^{\beta}$ and  
$\Omega_W=\oplus_{\beta\in L+\lambda_W } \Omega_W^{\beta}$ 
with $\Omega_W^{\beta}=W^\beta\cap \Omega_W$.
For $a\in \Omega_V^\alpha$, set 
$$
  Y_{\Omega_W}(a,z):=E^-(\alpha,z) Y_W(a,z) E^+(\alpha,z) z^{-\alpha_{(0)}}.
$$
Then it is shown in \cite{DL} and \cite{L4} that 
$(\Omega_W,Y_{\Omega_W}(\cd,z))$ is an irreducible $(L+\lambda_W)$-graded 
$\Omega_V$-module.
Moreover, we have the following theorem:

\begin{thm}\label{L4}
  (\cite[Theorem 3.16]{L4})
  For a rational VOA $V$, the associated $L$-graded generalized 
  vertex operator algebra $(\Omega_V,\w_\Omega)$ is also rational in the 
  sense that every $\Omega_V$-module with an $L$-set grading is completely 
  reducible.
  Moreover, the map which associates to an irreducible $V$-module $W$ 
  an irreducible $\Omega_V$-module $\Omega_W$ defines a bijection between
  the set of inequivalent irreducible $V$-modules and the set of inequivalent
  irreducible $\Omega_V$-modules with $L$-set gradings.
\end{thm}

\begin{rem}
  Even if $V$ is a simple VOA, the generalized VOA $\Omega_V$ may contain 
  a non-trivial ideal.
  However, the theorem above says that there is no {\it $L$-graded} ideal
  in $\Omega_V$.
\end{rem}

By this theorem, we can expect that the space of $q$-characters 
$\ch_{\Omega_W}(\tau)$, where $W$ runs over irreducible $V$-modules, has a 
modular invariance property.
Below we show modular transformation laws of $\ch_{\Omega_W}(\tau)$.

Let $a\in V^\alpha$ and $u,v\in H$.
Then $a_{(n)} W^\beta \subset W^{\alpha+\beta}$ for $\beta \in L+\lambda_W$, 
$n\in \Z$, whereas $u_{(0)}$ and $v_{(0)}$ acts on $W^{\beta}$ semisimply.
So we have $Z_W(a;(u,v);\tau)=0$ unless $\alpha =0$.
That is, our theta function is effective only for elements in 
$V^0=M_\h(1,0)\tensor \Omega_V^0$.

\begin{prop}
  Let $W$ be an untwisted $V$-module and $a\in V^\alpha$.
  If $\alpha\ne 0$, then $Z_W(a;(u,v);\tau)=0$.
\end{prop}

If $a\in \Omega_V^0\subset V^0$, then $Y_{\Omega_W}(a,z)=Y_W(a,z)$ 
so that $a$ acts on each $\Omega_W^\beta$.
Moreover, both $u_{(0)}$ and $v_{(0)}$ are also commutative with 
actions of $M_\h(1,0)$ on $W$. Therefore, we have 
$$
\begin{array}{l}
  Z_W(a;(u,v);\tau)
  = e^{\pii \la u,v\ra} \tr_W o(a) \sigma (v) q^{L(0)+u_{(0)}+\hf \la u,u\ra
  -\fr{1}{24}c}
  \vsb\\
  = e^{\pii \la u,v\ra} \tr_{M_\h(1,0)}q^{L_\h(0)-\dim \h /24} 
  \cd  \tr_{\Omega_W} o(a) \sigma (v) q^{L(0)-L_\h(0)+u_{(0)} 
  +\hf \la u,u \ra -\fr{1}{24}(c-\dim \h)}
  \vsb\\
  = e^{\pii \la u,v\ra} \eta(\tau)^{-\dim \h} 
  \cd \tr_{\Omega_W} o(a) \sigma (v) q^{L_\Omega (0)+u_{(0)} 
  +\hf \la u,u \ra -\fr{1}{24}(c-\dim \h)}, 
\end{array}
$$
where $L_\h(0)=(\w_\h)_{(1)}$ and $L_\Omega(0)=(\w -\w_\h)_{(1)}$ are
degree operators on $M_\h(1,0)$ and $\Omega_W$, respectively.
By the above equality, the essential ingredient of a theta function
$Z_W(a;(u,v);\tau)$ comes from the structure of the space 
$\Omega_W=\oplus_{\beta \in L+\lambda_W} \Omega_W^\beta$ 
of highest weight vectors.
Now for $a\in \Omega_V$ and $u,v\in H$, set a trace form on $\Omega_W$ by
\begin{equation}
  X_{\Omega_W}(a;(u,v);\tau)
  :=e^{\pii \la u,v\ra} \tr_{\Omega_W} o(a) \sigma (v) 
     q^{L_\Omega (0)+u_{(0)} +\hf \la u,u \ra -\fr{1}{24}(c-\dim \h)}.
\end{equation}
Then $Z_W(a;(u,v);\tau)= \eta(\tau)^{-\dim \h} \cd 
X_{\Omega_W}(a;(u,v);\tau)$.
Let $\{ W^1,\dots,W^N\}$ be the set of all inequivalent irreducible untwisted 
$V$-modules.
Then by Theorem \ref{main thm} we have the following modular transformation 
laws:

\begin{thm}
  For $a\in \Omega_V^0$ with $L_\Omega(0)a= \wt_\Omega(a)\cd a$, 
  $\wt_{\Omega}(a) \in \Q$, we have
  $$
  \begin{array}{l}
     X_{\Omega_{W^i}}(a;(u,v);\tau+1)
     = e^{\pii /12} \dsum_{j=1}^N T_{ij}\, X_{W^j}(a;(u,u+v),\tau),
     \\
     X_{\Omega_{W^i}}(a;(u,v);-1/\tau)
     = (-\sqrt{-1}\, \tau)^{\dim \h /2} \cd \tau^{\wt_\Omega(a)}\,
       \dsum_{j=1}^N S_{ij}\, X_{W^j}(a;(v,-u);\tau),
  \end{array}
  $$
  where $T_{ij}=A_{ij}(\binom{1\ 1}{0\ 1})$ and 
  $S_{ij}=A_{ij}(\binom{0\ -1}{1\ \ 0})$ are constants given by Theorem 4.2 
  as in the equation \eqref{transformation}.
\end{thm}

\pf
The following transformation laws are well-known:
$\eta(\tau+1)=e^{\pii /12}\eta(\tau)$ and 
$\eta(-1/\tau)=(-\sqrt{-1}\, \tau)^{1/2} \eta (\tau)$.
Since $L_\h(0) a =0$ by definition of $\Omega_V$, we have $L(0) a
= (L_\h(0)+L_\Omega (0))a= \wt_\Omega (a)\cd a$.
Thus by combining Theorem \ref{main thm} with the relation
$Z_{W^i}(a;(u,v);\tau)=\eta(\tau)^{-\dim \h} X_{\Omega_{W^i}}(a;(u,v);\tau)$, 
we have the desired equalities.
\qed

Note that $X_{\Omega_{W^i}}(\vac;(0,0);\tau) 
=\tr_{\Omega_{W^i}}q^{L_\Omega(0)-(c-\dim \h)/24}
=\ch_{\Omega_{W^i}}(\tau)$ is the $q$-character of an 
$\Omega_V$-module $\Omega_{W^i}$.
By the theorem above, the space of $q$-characters of $\Omega_V$-modules
is not invariant under $\SL_2(\Z)$ as we have an extra term $(-\sqrt{-1}\, 
\tau )^{-\dim \h /2}$.
However, we can eliminate the term $(-\sqrt{-1}\, \tau)^{-\dim \h /2}$ 
in the following way.
Assume that the rational lattice $L$ contains a positive definite even 
lattice $K$ such that $\mathrm{rank}(L)=\mathrm{rank}(K)(=\dim\h)$.
Then $V$ contains a lattice VOA $V_K$ associated to $K$, and 
all $W^i$ are twisted modules for $V_K$.
We further assume that we can choose a lattice $K$ such that all $W^i$ are 
{\it untwisted} $V_K$-modules.
How to choose such a lattice $K$ is shown in \cite{L4}, and such $K$ 
always satisfies $L+\lambda_{W^i} \subset K^\circ$ for all $1\leq i\leq 
N$.
Since $V_K$ is rational, we have the following decomposition for each 
$W^i$:
$$
  W^i=\bigoplus_{\mu +K\in (L+\lambda_{W^i})/K} V_{K+\mu}\tensor 
  \hom_{V_K} (V_{K+\mu},W^i).
$$
Set $\Omega_{W^i}^{\mu}:= \hom_{V_K}(V_{K+\mu},W^i)$ and 
$\Omega_{W^i}^K:=\oplus_{\mu +K \in (L+\lambda_{W^i})/K} \Omega_{W^i}^\mu$.
Then it is shown in \cite{DL} and \cite{L4} that there is an ideal $I$ of
$\Omega_V$ such that the quotient generalized vertex operator algebra
$\Omega_V/I$ is simple and isomorphic to $\Omega_V^K$.
So the space $\Omega_V^K$ naturally possesses a structure of a 
simple $L/K$-graded generalized vertex operator algebra.
Moreover, it is also shown in \cite{L4} that each $\Omega_{W^i}^K$ is 
an irreducible $(L+\lambda_{W^i})/K$-graded $\Omega_V^K$-module, and 
every $\Omega_V^K$-module with an $L/K$-set grading is a direct
sum of copies of $\Omega_{W^i}^K$'s.

On the other hand, we have 
$$
\begin{array}{ll}
  Z_{W^i}(\vac;(0,0);\tau)
  &=\sum_{\mu+K \in (L+\lambda_{W^i})/K} 
  \ch_{V_{K+\mu}}(\tau ) \cd \ch_{\Omega_{W^i}^\mu}(\tau)
  \vsb\\
  &= \sum_{\mu+K \in (L+\lambda_{W^i})/K} \eta (\tau)^{-\mathrm{rank}(K)}\cd 
    \theta_{K+\mu}(\tau)\cd \ch_{\Omega_{W^i}^\mu}(\tau)
  \vsb\\
  &= \eta (\tau)^{-\dim \h} \sum_{\mu+K \in (L+\lambda_{W^i})/K} 
    \theta_{K+\mu}(\tau) \cd \ch_{\Omega_{W^i}^\mu}(\tau),
\end{array}
$$
where $\theta_{K+\mu}(\tau)=\sum_{\beta \in K+\mu} q_\tau^{\la \beta,\beta 
\ra /2}$ are theta functions defined on the lattice $K$.
Thus 
$$
  X_{\Omega_{W^i}}(\vac;(0,0);\tau)=\sum_{\mu +K \in (L+\lambda_{W^i})/K} 
  \theta_{K+\mu}(\tau)\cd \ch_{\Omega_{W^i}^\mu}(\tau).
$$
Since the term $(-\sqrt{-1}\, \tau)^{-\dim \h/2}=(-\sqrt{-1}\, \tau
)^{-\mathrm{rank}(K)/2}$ also appears in the modular transformation of 
$\theta_{K+\mu}(\tau)$, we have a desired cancellation.
By this observation, it is very likely to happen that the space of 
$q$-characters 
$$
  \ch_{\Omega_{W^i}^K}(\tau)
  =\sum_{\mu +K \in (L+\lambda_{W^i})/K} \ch_{\Omega_{W^i}^\mu} (\tau)
$$
of $\Omega_V^K$-modules $\Omega_{W^i}^K$, $1\leq i\leq N$, are invariant
under the action of $\SL_2(\Z)$.
This question will be discussed in another paper.

\begin{center}
  \gs{Acknowledgment}
\end{center}
\begin{quote}
  The author would like to thank Professor Masahiko Miyamoto for 
  telling him the proof of Lemma 2.4 of \cite{M3}.
\end{quote}

\small
\setlength{\baselineskip}{12pt}

\end{document}